\documentclass[10pt,reqno]{amsart}

\newcommand{\myreferences}{references}

\usepackage[utf8]{inputenc}
\usepackage{amsmath}
\usepackage{amssymb}
\usepackage{hyperref}
\usepackage{stackrel}
\usepackage{xcolor}
\usepackage{graphicx}
\usepackage{subfigure}
\usepackage{comment}
\usepackage{enumerate}

\newtheorem{theorem}{Theorem}[section]

\newtheorem{definition}[theorem]{Definition}

\newtheorem{proposition}[theorem]{Proposition}
\newtheorem{remark}[theorem]{Remark}
\newtheorem{corollary}[theorem]{Corollary}

\newtheorem{correction}[theorem]{Corrections}

\DeclareMathOperator{\spann}{span}

\newcommand{\atlas}{\mathcal{A}}

\newcommand{\geodesicspray}{\mathbf{G}}
\newcommand{\geodesicspraycotangent}{\mathcal{G}}

\newcommand{\familiavetoresunitarios}{\mathcal{F}}

\newcommand{\slitcotangentbundle}{T^{\ast}M\setminus 0}
\newcommand{\slitcotangentbundleU}{T^{\ast}U\setminus 0}
\newcommand{\slittangentbundle}{TM\backslash 0}

\newcommand{\LegendreTransform}{\mathcal{L}}

\makeatletter
\newcommand{\tpitchfork}{%
  \vbox{
    \baselineskip\z@skip
    \lineskip-.52ex
    \lineskiplimit\maxdimen
    \m@th
    \ialign{##\crcr\hidewidth\smash{$-$}\hidewidth\crcr$\pitchfork$\crcr}
  }%
}
\makeatother

\makeatletter
\@namedef{subjclassname@2020}{\textup{2020} Mathematics Subject Classification}
\makeatother

\title[Addendum to ``Geodesic fields for Pontryagin type ... '']{Addendum to ``Geodesic fields for Pontryagin type $C^0$-Finsler manifolds''}

\author{Hugo Murilo Rodrigues and Ryuichi Fukuoka}

\subjclass[2020]{Primary 49J15, 53B40 Secondary 53C22}

\begin{document}

\maketitle

\begin{abstract}
In the paper ``Geodesic fields for Pontryagin type $C^0$-Finsler manifolds'', the same authors introduced a geometric object called Pontryagin-type $C^0$-Finsler structures, which includes Finsler structures on a smooth manifold $M$.
Such a structure is defined locally by a control system of $C^1$ unit vector fields on an open subset $U$ of $M$, such that Pontryagin's maximum principle (PMP) is applicable. 
As a result of PMP, we obtain the extended geodesic field $\mathcal{E}$, which is a multivalued vector field on the slit cotangent bundle $\slitcotangentbundleU = \{(x,\alpha): x\in U, \alpha \in T^{\ast}_x U\}$, such that every minimizing path $x(t)$ on $(U,F)$ is an extremal (projection of a Pontryagin extremal $(x(t),\alpha(t))$).    
Some time after its introduction, its use and external comments, it is useful to provide remarks that clarify the applicability and limitations of this theory.
In addition, we point out two errors contained in that work which, although discernible from the context, we believe should be pointed out.
Finally, we prove some cases where the property that a minimizing path is the projection of a Pontryagin extremal can be extended beyond a coordinate system, sometimes to $M$.  
We also work on the problem of determining when $\mathcal{E}$ does not depend on the control system chosen to represent $F$. 
In particular, $\mathcal{E}$ satisfies these two properties globally for left-invariant $C^0$-Finsler structures on Lie groups, and for $C^1$-Finsler structures (without the condition of strong convexity).
\end{abstract}

\section{Introduction}

The aim of this addendum to manuscript \cite{Fukuoka-Rodrigues} is to provide remarks that clarify the applicability and limitations of its results. 
We also point out two errors contained in this work which, although discernible from the context, we believe should be emphasized.
We also introduce a result that extends the applicability of the theory developed in \cite{Fukuoka-Rodrigues} to two classes of Pontryagin-type $C^0$-Finsler manifolds.

As usual, $M$ is a smooth manifold, $T_xM$ and $T^{\ast}_xM$ are the tangent and cotangent spaces of $M$ at $x \in M$ respectively, $TM$ and $T^{\ast}M$ are the tangent and cotangent bundle of $M$.
A $C^0$-Finsler structure on $M$ is a continuous function $F:TM \rightarrow \mathbb{R}$ such that $F_x = F(x,\cdot):T_xM \rightarrow \mathbb{R}$ is an asymmetric norm for every $x \in M$.
$F_{\ast} : T^\ast M \rightarrow \mathbb{R}$ is the fiberwise dual $C^0$-Finsler structure.
The unit spheres, open balls and closed balls in $(T_xM,F)$, centered at $a$ and radius $r$, are denoted by $S_F[x,a,r]$, $B_F(x,a,r)$ and $B_F[x,a,r]$ respectively. 
When $a=0$ and $r=1$, we denote them simply by $S_F[x]$, $B_F(x)$ and $B_F[x]$ respectively. 
Analogous notations for spheres and balls are used in similar cases. 

The geodesics of a (smooth) Finsler manifold $(M,F)$ can be studied through its geodesic spray $\geodesicspray$, which is a smooth vector field on $\slittangentbundle$.
We can also use the Legendre Transform $\mathcal{L}:\slittangentbundle \rightarrow \slitcotangentbundle$ to represent the geodesic spray on $\slitcotangentbundle$ as $\geodesicspraycotangent = d\LegendreTransform (\geodesicspray)$.
Since $\LegendreTransform$ is a diffeomorphism in Finsler geometry, the structures $\geodesicspray$ and $\geodesicspraycotangent$ are equivalent (see \cite{BaoChernShen}).

In \cite{Fukuoka-Rodrigues}, the authors studied to what extent we can generalize $\geodesicspray$ and/or $\geodesicspraycotangent$ for $C^0$-Finsler manifolds. 
Since variational tools aren't available for general $C^0$-Finsler structures, the authors considered a certain type of horizontal smoothness, where the Pontryagin's maximum principle (PMP) is applicable (see Definition \ref{definition_Pontryagin_type}).

In this work, we also introduce formal definitions for certain previously given but unnamed concepts, such as the $F$-family of controlled $C^1$ unit vector fields, due to their frequent use and importance in this work. 

\begin{definition}
\label{definition_Pontryagin_type}
A $C^0$-Finsler manifold is of Pontryagin type at $p \in (M,F)$ if there exist a neighborhood $U$ of $p$, a coordinate system $\phi = (x^1, \ldots, x^n): U \rightarrow \mathbb R^n$ (with the respective natural coordinate system $\phi_{TU} := d\phi: (x^1, \ldots, x^n,$ $ y^1, \ldots, y^n): TU \rightarrow \mathbb R^{2n}$ of the tangent bundle) and a family of $C^1$ unit vector fields 
\[
\familiavetoresunitarios_U = \{x \mapsto X_u(x) = (y^1(x^1, \ldots, x^n,u), \ldots, y^n(x^1, \ldots, x^n,u)); u \in S^{n-1}\}
\]
on $U$ parameterized by $u \in S^{n-1}$ such that
\begin{enumerate}
\item $u \mapsto X_u(x)$ is a homeomorphism from $S^{n-1} \subset \mathbb R^n$ onto $S_F[x,0,1]$ for every $x \in U$;
\item $(x,u) \mapsto (y^1(x^1, \ldots, x^n,u), \ldots, y^n(x^1, \ldots, x^n,u))$ is continuous;
\item $(x,u) \mapsto (\frac{\partial y^1}{\partial x^i} (x^1, \ldots, x^n, u), \ldots, \frac{\partial y^n}{\partial x^i} (x^1, \ldots, x^n, u))$ is continuous for every $i=1, \ldots, n$.
\end{enumerate} 
$\mathcal{F}_U$ is an {\em $F$-family of controlled $C^1$ unit vector fields defined on $U$ (or simply an $F$-family on $U$)}.
We say that $F$ is of Pontryagin type on $M$ if it is of Pontryagin type at every $p \in M$. 
In this case $(M,F)$ is a Pontryagin-type $C^0$-Finsler manifold.
\end{definition}

The family $\familiavetoresunitarios_U$ determines a control system that meet the minimum requirements for applying PMP. 
Moreover, the definition of Pontryagin type doesn't depend on the choice of the coordinate system (see \cite{Fukuoka-Rodrigues}). 

The study of minimizing paths of $(U,F)$ through $\familiavetoresunitarios_U$ is a time minimizing problem, because the elements of $\familiavetoresunitarios_U$ are unit vector fields. 
PMP is applied here as follows: 
For every control $u\in S^{n-1}$, define the Hamiltonian $H_u (x,\alpha) = \alpha (X_u(x))$ with respect to $u$ on $\slitcotangentbundleU = \{(x,\alpha): x\in U, \alpha\in T^{\ast}_x M \backslash \{0\}\}$, and the Hamiltonian vector field $\vec{H}_u$ on $\slitcotangentbundleU$ through the equation $d H_u = \omega (\cdot,\vec{H}_u)$, where $\omega$ is the canonical symplectic 2-form on $T^{\ast}U$.
Now we are in position to define the extended geodesic field, which is a consequence of PMP:

\begin{definition}
\label{definition_extended_geodesic_field}
Let $(M,F)$ be a Pontryagin-type $C^0$-Finsler manifold and $\mathcal{F}_U$ an $F$-family of $(U,F)$. 
The extended geodesic field of $(U,F)$ is the rule $\mathcal E$ that associates each $(x,\alpha) \in T^*U\backslash 0$ to the set 
\begin{align}
\label{align_extended_gedesic_field}
\mathcal E(x, \alpha) = \{\vec H_u(x, \alpha);u \in \mathcal C(x,\alpha)\},
\end{align} 
where $\mathcal C(x,\alpha) = \{u \in S^{n-1};$ $ \alpha(X_u(x)) = \max_{v\in S^{n-1}} \alpha(X_v(x)) \}$.  Pontryagin extremals are absolutely continuous curves $\gamma:[a,b] \rightarrow T^\ast U\backslash 0$ which are solutions of the differential inclusion
\begin{align}
\label{align_differential_inclusion}
\gamma^\prime(t) \in \mathcal E(\gamma(t)).
\end{align}
Extremals $x(t)$ are projections of Pontryagin extremals $\gamma(t) = (x(t),\alpha(t))$ on $U$.
\end{definition}

Let us write $\vec{H}_u(x,\alpha)$ (and, consequently, $\mathcal{E}(x,\alpha)$) in a coordinate system. 
Let $(x^1, \ldots, x^n, \alpha_1, \ldots, \alpha_n)$ be the coordinate system on $T^\ast U$ induced by $\phi$. 
Henceforth, the convention of summing over upper and lower indices applies.
If we write 
\begin{align}
  X_{u}(x) = f^i(x,u) \frac{\partial}{\partial x^i}, \nonumber
\end{align}
the Hamiltonian $H_u$ with respect to $u \in S^{n-1}$ is given by
\begin{align}
  H_{u}(x,\alpha) = \alpha_i f^i(x,u),
\end{align}
and the corresponding Hamiltonian vector field equals
\begin{align}
\vec{H}_u(x,\alpha) 
& = f^i(x, u) \frac{\partial}{\partial x^i} - \frac{\partial}{\partial x^i}\left( H_u (x,\alpha) \right) \frac{\partial}{\partial \alpha_i} \nonumber \\
& = f^i(x, u) \frac{\partial}{\partial x^i} - \alpha_j \frac{\partial f_j}{\partial x^i} (x,u) \frac{\partial}{\partial \alpha_i}. \label{align_Hamiltonian_vector_field_coordinates}
\end{align}

The next theorem states that every minimizing path parametrized by arc length is a controlled path with respect to the control system $\mathcal{F}_U$.
Therefore, minimizing paths of $(U,F)$ are among its extremals (see \cite{Fukuoka-Rodrigues}).

\begin{theorem}
  \label{theorem_minimizing} Let \((M, F)\) be a Pontryagin-type \(C^0\)-Finsler manifold. 
  Then every minimizing curve \(x(t)\) of \((U, F)\) parameterized by arclength is the projection of a \textit{Pontryagin extremal} \((x(t), \alpha(t))\). 
  Moreover, the \textit{Hamiltonian} \(F_{\ast} (x(t),\alpha(t)) = H_{u(t)}(x(t), \alpha(t))\) is constant.
\end{theorem}

Therefore, if we want to study minimizing paths $x(t)$ of $(M,F)$ with respect to a coordinate system $\phi: U \rightarrow \mathbb{R}^n$ in a neighborhood of $p \in M$ as it is usually done in Finsler geometry, we can fix an $F$-family $\mathcal{F}_U$ on $U$ (eventually choosing a smaller neighborhood $U$), and consider the differential inclusion $(\dot{x}(t),\dot{\alpha}(t)) \in \mathcal{E}(x(t),\alpha(t))$ on $T^\ast U \backslash 0$, which in the coordinate system is given by \eqref{align_extended_gedesic_field}, \eqref{align_differential_inclusion} and \eqref{align_Hamiltonian_vector_field_coordinates}.

\begin{correction}
Here we would like to make corrections in Definition 4.2 and Theorem 4.4 of \cite{Fukuoka-Rodrigues}. 
For Definition 4.2, it only makes sense to define $\mathcal{E}$ on $T^{\ast}U \backslash 0$, where $U$ is a coordinate open subset where the $F$-family is defined. 
In the same way, Theorem 4.4 makes sense only for minimizing paths parametrized by arc length defined on $U$.
\end{correction}

If we are interested in extending Theorem \ref{theorem_minimizing} to domains of $M$ that are larger than coordinate open subsets where $F$-families are defined, then the following remark is useful and straightforward (see, e.g. \cite{Prudencio-Fukuoka-1} and \cite{Fukuoka-Setti-2}).   

\begin{remark}
  \label{remark_importance_of_the_domain_of_a_F_family}
  Notice that $\mathcal{E}$ doesn't depend on the choice of the coordinate system.
  Therefore, if we have a family $\familiavetoresunitarios_V= \{X_u: V \rightarrow TV, u\in S^{n-1}\}$ of controlled $C^1$ unit vector fields defined on an open subset $V$ of $M$ ($V$ not being necessarily a coordinate open subset of $M$), such that its restriction $\familiavetoresunitarios_U := \{X_u\vert_U, u \in S^{n-1}\}$ to every coordinate open subset $U$ of $V$ satisfies the conditions of Definition \ref{definition_extended_geodesic_field}, then $\mathcal{E}$ can be defined on $T^\ast V \backslash 0$. 
  In addition, the proof of Theorem 4.4 of \cite{Fukuoka-Rodrigues}, that every minimizing path parametrized by length is a controlled path, 
  works without changes on $V$. 
  Therefore, every minimizing path on $V$ is an extremal (see \cite{Sachkov}).
  Such a family $\mathcal{F}_V$ is also called an {\em $F$-family on $V$}.
  For instance, in \cite{Prudencio-Fukuoka-1} and \cite{Fukuoka-Setti-2}, the authors considered the $F$-family of unit left-invariant vector fields on Lie groups $G$, which are defined globally on $G$.
\end{remark}

The $F$-families on open subsets $V \subset M$ are part of the definition of Pontryagin-type $C^0$-Finsler structures, and $\mathcal{E}$ depends, a priori, on the chosen $F$-family.
In some classical cases, as in Finsler manifolds $(M,F)$, $\mathcal{E}$ doesn't depend on the $F$-families because 
$\mathcal{E}=d\mathcal{L}(\mathbf{G})$ (see \cite{Fukuoka-Rodrigues}).
For a fixed $(x,\alpha)$, if $\mathcal{E}(x,\alpha)$ doesn't depend on the $F$-family, we say that {\em $\mathcal{E}(x,\alpha)$ depends only on $F$}.
If this property holds for every $(x,\alpha) \in T^{\ast}V\backslash 0$, we say that $\mathcal{E}$ {\em depends only on $F$ on $T^{\ast}V\backslash 0$}.
If $V=M$, we simply say that {\em $\mathcal{E}$ depends only on $F$}.
In Section \ref{section_main_theorem}, we prove a sufficient condition for $\mathcal{E}(x,\alpha)$ to depend only on $F$.

In Section \ref{section_examples}, we prove that $\mathcal{E}$ depends only on $F$ for $G$-invariant $C^0$-Finsler structures on homogeneous spaces, as well as for Pontryagin-type $C^0$-Finsler structures such that $F_{x}: T_{x}M \backslash \{0\} \rightarrow \mathbb{R}$ is a $C^1$ map for every $x \in M$.
Another point that is addressed is to what extent minimizing paths in $M$ are projections of Pontryagin extremals (cf. Theorem 4.4 of \cite{Fukuoka-Rodrigues}).
We prove that this property holds for Pontryagin-type $C^0$-Finsler structures such that $F_x: T_xM\backslash \{0\} \rightarrow \mathbb R$ is a $C^1$ map for every $x \in M$.

\begin{remark}
  Generative AI tools were used only for text readability improvements.
  The mathematical content was exclusively developed by the authors.
\end{remark}

\section{Main theorem}

\label{section_main_theorem}

This section is devoted to the proof of a sufficient condition for $\mathcal{E}(x,\alpha)$ to depend only on $F$.
For a coordinate system $\phi = (x^1, \ldots, x^n)$ on $U \subset M$, consider the natural coordinate systems $(x^1, \ldots, x^n, y^1, \ldots, y^n)$ and $(x^1, \ldots, x^n, \alpha_1, \ldots, \alpha_n)$ on $TU$ and $T^\ast U$ respectively. 

\begin{theorem}
  \label{theorem_main}
  Let $(M,F)$ be a Pontryagin-type $C^0$-Finsler manifold and consider $(x_0,\alpha_0) \in \slitcotangentbundle$. 
  Let $\mathcal{F}_U = \{X_u: U \rightarrow TU; u \in S^{n-1}\}$ be an $F$-family defined on a neighborhood $U$ of $x_0$, and for every $v \in S_F[x_0]$, denote by $u(v) \in S^{n-1}$, the control such that $X_{u(v)}(x_0)=v$.
  Suppose that for every $v_0 \in \mathcal{C}(\alpha_0)$, there exist
  a $1$-form $\beta: U \mapsto T^\ast U \setminus 0$, continuous at $x_0$,  
  such that:
\begin{enumerate}
  \item $\beta(x_0) = \alpha_0$; \label{item_alpha_x_0}
  \item \label{item_alpha_X} The derivative of $x \in U \mapsto \beta(X_{u(v_0)}(x))$ at $x_0$ is the zero functional; 
  \item \label{item_F_star_alpha} The derivative of $x \in U \mapsto F_{\ast}(x,\beta(x))$ at $x_0$ is the zero functional.
\end{enumerate}  
Then $\mathcal{E}(x_0,\alpha_0)$ depends only on $F$.
\end{theorem}

\begin{remark} 
  Items (\ref{item_alpha_X}) and (\ref{item_F_star_alpha}) of Theorem \ref{theorem_main} are quite technical.
  Their use become more concrete in the examples of Section \ref{section_examples}.
\end{remark}

\begin{proof}[Proof of Theorem \ref{theorem_main}]
Consider a coordinate system, given by a $C^1$-diffeomorphism $\phi = (x^1, \ldots,$ $x^n):U \rightarrow \mathbb{R}^n$ defined in a neighborhood $U$ of $x_0$ such that 
\begin{align}
  \label{align_Xuv0}
  X_{u(v_0)} = \frac{\partial}{\partial x^1} \text{ on }U
\end{align} 
and $\ker \alpha_0 = \spann \left\{\frac{\partial}{\partial x^2}, \ldots, \frac{\partial}{\partial x^n} \right\}$ (Here we can choose a smaller neighborhood $U$ of $x_0$ if it is necessary).
In this coordinate system, the Hamiltonian is given by $H_{u(v_0)}(x,\alpha)=\alpha(X_{u(v_0)}) = \alpha_1 > 0$ and 
\begin{align}
\vec{H}_{u(v_0)} (x, \alpha) 
= f^i(x,u(v_0))\frac{\partial}{\partial x^i} - \alpha_j \frac{\partial f^j}{\partial x^i}(x,u(v_0)) \frac{\partial}{\partial \alpha_i} 
= \frac{\partial}{\partial x^1}.
\label{align_Hamiltonian_vector_field_u1}
\end{align}
Consider $u_2, \ldots, u_n \in S^{n-1}$ such that $X_{u_j}(x_0)= \frac{\partial}{\partial x^j}(x_0)$, $j=2,\ldots, n$. 
Set $u_1=u(v_0)$.
If necessary, consider a smaller $U$ such that $\{X_{u_1}(x), \ldots, X_{u_n}(x)\}$ is a basis of $T_xU$ for every $x \in U$.

Now consider another $F$-family $\tilde{\mathcal{F}}_U = \{X_{w}: U \rightarrow T U; w \in S^{n-1}\}$ (We can always suppose that these families of vector fields are defined in the same neighborhood $U$ of $x_0$, taking the intersection of the domains whenever necessary).
Consider the control $w(v_0) \in S^{n-1}$ such that $X_{w(v_0)}(x_0)  = \frac{\partial}{\partial x^1}$.
Write 
\[
X_{w(v_0)}(x) = a^i(x) X_{u_i}(x),
\]
and notice that 
\begin{align}
  \label{align_ai}
  a^1(x_0)=1 \text{ and } a^i(x_0)=0 \text{ for } i=2,\ldots, n.
\end{align}
Define the $1$-form $\tilde{\alpha}_0 = (\alpha_0)_1 dx^1$ on $U$.
Then 
\[
H_{w(v_0)}(x, (\alpha_{0})_1,0,\ldots, 0) = (\tilde{\alpha}_0 (X_{w(v_0)}))(x) 
= \tilde{\alpha}_0 (a^i(x) X_{u_i}(x)) 
= a^i(x) \tilde{\alpha}_0(X_{u_i}(x)).
\]
By \eqref{align_Hamiltonian_vector_field_coordinates}, the vertical part of $\vec{H}_{w(v_0)}$ is given by
\begin{align}
\label{align_vertical_Hamiltonian_w1_alpha0}
& - \frac{\partial}{\partial x^j}\left(  H_{w(v_0)} (x,(\alpha_0)_1, 0, \ldots, 0) \right) \nonumber \\
= & -\frac{\partial a^i}{\partial x^j} (x) . \tilde{\alpha}_0(X_{u_i}(x))-a^i(x)\, \frac{\partial}{\partial x^j}\left( \tilde{\alpha}_0(X_{u_i}(x)) \right).
\end{align}
Therefore, $\vec{H}_{u(v_0)}(x_0, \alpha_0) = \vec{H}_{w(v_0)}(x_0, \alpha_0)$ iff 
\begin{align}
\frac{\partial a^1}{\partial x^j}(x_0) = 0
\label{align_vertical_nulo}
\end{align}
for every $j=1,\ldots, n$, due to \eqref{align_Xuv0} and \eqref{align_ai}. 
Let us rewrite \eqref{align_vertical_nulo}. Set
\[
X_{u_i}(x)=X^j_{u_i}(x)\frac{\partial}{\partial x^j}.
\]
Then
\[
X_{w(v_0)}(x):= f^j (x,w(v_0))\frac{\partial}{\partial x^j} = a^i(x) X^j_{u_i}(x)\frac{\partial}{\partial x^j},
\]
what implies
\[
\frac{\partial f^j}{ \partial x^k}(x,w(v_0)) = \frac{\partial a^i}{\partial x^k}(x) . X^j_{u_i}(x) + a^i(x).\frac{\partial X^j_{u_i}}{\partial x^k}(x)
\]
and
\[
\frac{\partial f^j}{ \partial x^k}(x_0,w(v_0)) = \frac{\partial a^j}{\partial x^k}(x_0).
\]
Therefore, $\vec{H}_{u(v_0)}(x_0, \alpha_0) = \vec{H}_{w(v_0)}(x_0, \alpha_0)$ iff 
\begin{align}
\frac{\partial f^1}{\partial x^k}(x_0, w(v_0)) = 0
\label{align_vertical_nulo_2}
\end{align}
for every $k=1,\ldots, n$, due to \eqref{align_vertical_nulo}.

Write $\beta(x) = \beta_i(x) dx^i$, where $\beta_1(x_0)=1$, and $\beta_i(x_0)=0$ for $i=2,\ldots, n$.
Notice that the map defined in Item \eqref{item_alpha_X} is given by $x \mapsto \beta_1(x)$, and
$\beta_1(x) = \beta_1(x_0) + o(\Vert x-x_0\Vert)$ when $x \rightarrow x_0$ due to its derivative at $x_0$.
Suppose that there exists $\tilde{v} = \tilde{v}^s\frac{\partial}{\partial x^s} \neq 0$ such that the directional derivative of $f^1$ satisfies $\frac{\partial f^1}{\partial x^s}(x_0,w_1) \tilde{v}^s > 0$.
We claim for a contradiction.
Then Item \eqref{item_F_star_alpha} implies
\begin{align}
  0 & = \lim_{t \rightarrow 0} \frac{F_\ast (x_0 + t\tilde{v}, \beta (x_0 + t\tilde{v})) - F_{\ast}(x_0,\beta(x_0))}{t } \nonumber \\ 
  & = \lim_{t\rightarrow 0} \frac{\sup_{u\in S^{n-1}} \beta(x_0 + t\tilde{v}) (X_u(x_0 + t\tilde{v})) - \beta_1(x_0)}{t} \nonumber \\
  & \geq \lim_{t\rightarrow 0} \frac{\beta(x_0 + t\tilde{v}) (X_{w(v_0)}(x_0 + t\tilde{v})) - \beta_1(x_0)}{t} \nonumber \\
  & = \lim_{t \rightarrow 0} \frac{\beta_i (x_0 + t\tilde{v}) f^i(x_0 + t\tilde{v},w(v_0))-\beta_1(x_0)}{t} 
  \nonumber \\
  & = \beta_1(x_0)\,\frac{\partial f^1}{\partial x^j}(x_0,w(v_0))\tilde{v}^j > 0 \nonumber
\end{align}
due to 
\begin{align}
  & \beta_1(x_0 + t\tilde{v}) f^1(x_0 + t\tilde{v},w(v_0)) + \sum_{k=2}^n \beta_k(x_0 + t\tilde{v}) f^k(x_0 + t\tilde{v},w(v_0)) \nonumber \\ 
  & = (\beta_1(x_0)+o(\vert t\vert)) \left(1+\frac{\partial f^1}{\partial x^j}(x_0,w(v_0)) \tilde{v}^j t + o(\vert t\vert)\right) \nonumber \\
  & + \sum_{k=2}^n o(\vert 1\vert)\left(\frac{\partial f^k}{\partial x^j}(x_0,w(v_0)) \tilde{v}^j t + o(\vert t\vert)\right), \nonumber
\end{align}
what gives a contradiction. 
Therefore, $\frac{\partial f^1}{\partial x^j}(x_0,w(v_0)) = 0$ for every $j=1,\ldots, n$, $\vec H_{u(v_0)}(x_0,\alpha_0) = \vec H_{w(v_0)}(x_0,\alpha_0)$, and $\mathcal{E}(x_0,\alpha_0)$ depends only on $F$.
\end{proof}

\section{To what extent $\mathcal{E}(x,\alpha)$ generalizes the geodesic spray}
\label{section_examples}

In this section, we present some particular cases of Pontryagin-type $C^0$-Finsler manifolds that satisfies Theorem \ref{theorem_main} for every $(x,\alpha) \in \slitcotangentbundle$, that is, $\mathcal{E}$ depends only on $F$.
In addition, we consider some cases where minimizing paths on $M$, parameterized by arc length, are extremals.

\subsection{Homogeneous manifolds}
This subsection follows Section 8 of \cite{Fukuoka-Rodrigues}, with minor adjustments to the notation.

Let $G$ be a Lie group, $H$ a closed subgroup of $G$, $a:G \times G/H \rightarrow G/H$ the canonical left-action and $\pi:G \rightarrow G/H$ the projection map. 
Denote $gH = \pi(g)$ and $a_g = a(g,\cdot)$. 
A $C^0$-Finsler structure $F$ on $G/H$ is $G$ invariant if $a_g$ is an isometry of $(G/H,F)$ for every $g \in G$.
In \cite{Fukuoka-Rodrigues}, it is proved that $G$-invariant $C^0$-Finsler structures are of Pontryagin type.

Let $\mathfrak{g}$ be the Lie algebra of $G$ and $\mathfrak{h}$ the Lie subalgebra corresponding to $H$.
Consider a subspace $\mathfrak{m}$ of $\mathfrak{g}$ such that $\mathfrak{g} = \mathfrak{m} \oplus \mathfrak{h}$.
Since $\pi \circ \exp \vert_{\mathfrak{m}}$ is a local diffeomorphism at the origin, fix a neighborhood $U$ of $0 \in \mathfrak{m}$ such that $\pi \circ \exp \vert_U: U \rightarrow G/H$ is a diffeomorphism over its image.
Set $\sigma = \pi \circ \exp \vert_{U}: U \rightarrow (\pi \circ \exp)(U)$ and $\sigma_g = a_g \circ \sigma : U \rightarrow a_g(\sigma(U))$.
Consider an identification of $\mathfrak{m}$ with $\mathbb{R}^n$.
It follows that the family
\begin{align}
  \atlas = \left\{ \sigma_g: U \rightarrow G/H: \sigma_g(x)=g \pi (\exp (x)) \right\}_{g\in G} \nonumber
\end{align}
is an atlas of $G/H$.

Consider $S = d(\sigma^{-1})_{eH}(S_F[eH])$ $\subset \mathfrak{m}$, where $e$ is the identity element of $G$.
We replace the control region $S^{n-1}$ by $S$ (see \cite{Fukuoka-Rodrigues}).
Then there exist a unique asymmetric norm $F_S$ on $\mathfrak{m}$ such that $S$ is the unit sphere of $(\mathfrak{m},F_S)$, and we have that $d\sigma_0: (\mathfrak{m},F_S) \rightarrow (T_{eH}(G/H),F_{eH})$ is an isometry.
Moreover, if $g \in G$, then $d(\sigma_g)_{0} = d(a_g)_{eH} \circ d\sigma_0 : (\mathfrak{m},F_S) \rightarrow (T_{gH}(G/H), F_{gH})$ is an isometry because $F$ is $G$-invariant. 

Fix $g\in G$. 
Although the notation will become longer, we choose to represent $\sigma_g(U)$ by $\{x=g \exp(v)H: v \in U\}$ and $S_F [g\exp(v)H]$ by $\{d(\sigma_{g\exp(v)})_{0}(u)\}_{u\in S}$.
The $F$-family on $\sigma_g(U)$ is defined by
\begin{align}
\mathcal{F}_{\sigma_g(U)} = \{ g \exp (v) H \in \sigma_{g}(U) \mapsto X_u(g\exp(v)H) = d(\sigma_{g\exp(v)})_{0}(u)\}_{u\in S}. \label{align_Ffamily_homogeneous}
\end{align}

\begin{proposition}[Homogeneous case]
  \label{proposition_homogeneous_case}
  Let $G$ a Lie group, $H$ a closed subgroup of $G$ and $G/H$ be the homogeneous manifold endowed with a $G$-invariant $C^0$-Finsler structure $F$. 
  Suppose that we are in the settings of this subsection. Then $\mathcal{E}$ depends only on $F$.
\end{proposition}

\begin{proof}
  Let $(g\exp(v_0)H,\alpha_0) \in T^\ast (G/H) \backslash 0$. 
  We will prove that $\mathcal{E}(g\exp(v_0)H,\alpha_0)$ depends only on $F$ using Theorem \ref{theorem_main}. 
  Consider $d(\sigma_{g\exp(v_0)})_{0}(u_0) \in \mathcal{C}(\alpha_0)$. 
  Let $d(\sigma_{g\exp (v_0)})_0^{\ast}:T^\ast_{g\exp(v_0)H}(G/H) \rightarrow \mathfrak{m}^\ast$ be the transpose of $d(\sigma_{g\exp (v_0)})_0$ and denote $\alpha_{\mathfrak{m}} = d(\sigma_{g \exp(v_0)})^\ast_0 (\alpha_0)$.
Consider the smooth $1$-form $\beta$ on $\sigma_g(U)$ given by 
\begin{align}
\beta(g\exp(v)H) = d((\sigma_{g\exp(v)})^{-1})^{\ast}_{g\exp(v)H}(\alpha_{\mathfrak{m}}). \label{align_1form_homogeneous}
\end{align}
  Item \eqref{item_alpha_x_0} follows from the definition of $\beta$.
  For Items \eqref{item_alpha_X} and \eqref{item_F_star_alpha}, notice that
  \begin{align}
    \beta(g\exp(v)H) (X_{u_0}(g\exp(v)H))  = \beta_\mathfrak{m}(u_0) \label{align_beta_calculated_on_Ffamily}
  \end{align}
  due to \eqref{align_Ffamily_homogeneous} and \eqref{align_1form_homogeneous}.
  Therefore, Item \eqref{item_alpha_X} holds because \eqref{align_beta_calculated_on_Ffamily} is constant.
  Moreover,
  \begin{align}
    F_\ast (g\exp(v)H,\beta(g\exp(v)H)) & = \sup_{u \in S} 
    \beta(g\exp(v)H)(d(\sigma_{g\exp(v)})_{0} (u)) \nonumber \\ & = \sup_{u \in S}\beta_{\mathfrak{m}}(u) = F_{S \ast}(\beta_{\mathfrak{m}}) \label{align_Fast_calculated_for_Item3}
  \end{align}
  for every $g\exp(v)H \in \sigma_{g}(U)$, and Item \eqref{item_F_star_alpha} holds because \eqref{align_Fast_calculated_for_Item3} is constant. 
\end{proof}

\begin{corollary}
  Let $G$ be a Lie group endowed with a left-invariant $C^0$-Finsler structure. 
  Then $\mathcal{E}$ depends only on $F$, and every minimizing path of $G$ parametrized by arc length is an extremal.
\end{corollary}

\begin{proof}
  This is a direct consequence of Proposition \ref{proposition_homogeneous_case} and the fact that the $F$-family of unit left-invariant vector fields are defined on $G$ (See Remark \ref{remark_importance_of_the_domain_of_a_F_family}).
\end{proof}

\subsection{Unique supporting hyperplane}

Fix $(x,\alpha) \in \slitcotangentbundle$. 
Notice that $v \in \mathcal{C}(\alpha)$ iff $\{\alpha = F_{\ast}(x,\alpha)\}$ is a supporting hyperplane of $B_{F}[x,0,1]$ at $v$.

\begin{proposition}[Unique supporting hyperplane]
  \label{proposition_unique_hyperplane}
  Let $(M,F)$ be a Pontryagin-type $C^0$-Finsler manifold and $(x_0,\alpha_0) \in \slitcotangentbundle$.
  Suppose that for every $v \in \mathcal{C}(\alpha_0)$, $B_F[x_0,0,1]$ admits a unique supporting hyperplane at $v$.
  Then $\mathcal{E}(x_0,\alpha_0)$ depends only on $F$.
\end{proposition}

\begin{proof}
  Let $\mathcal{F}_U = \{X_u:U \rightarrow TU; u\in S^{n-1}\}$ be an $F$-family defined on a neighborhood $U$ of $x_0$ and $v \in \mathcal{C}(\alpha_0)$.
  Consider a coordinate system, given by a $C^1$-diffeomorphism $\phi=(x^1, \ldots, x^n):U \rightarrow \mathbb{R}^n$, such that $X_{u(v)} = \frac{\partial}{\partial x^1}$ on $U$ and $\ker \alpha_0 = \spann \left\{ \frac{\partial}{\partial x^2}, \ldots, \frac{\partial}{\partial x^n} \right\}$.
  Consider $\beta:U \rightarrow T^\ast U$ be a $1$-form on $U$ such that $\{ \beta(x) = F_{\ast}(x_0,\alpha_0) \}$ is a supporting hyperplane of $B_F[x,0,1]$ at $X_{u(v)}(x)$ (this functional $\beta(x)$ isn't necessarily unique). 
  Notice that $\beta(x_0) = \alpha_0$ due to the uniqueness of the supporting hyperplane at $v$.
  Therefore, we have a (not necessarily continuous) $1$-form $\beta$, satisfying  Items \eqref{item_alpha_X} and \eqref{item_F_star_alpha} of Theorem \ref{theorem_main}.
  We claim that $\beta$ is continuous at $x_0$, what proves that $\mathcal{E}(x_0,\alpha_0)$ depends only on $F$ due to Theorem \ref{theorem_main}.
  
  We claim for a contradiction, supposing that $\beta$ isn't continuous at $x_0$.
  Consider the Euclidean (Riemannian) metric $\left< \cdot, \cdot \right>$ on $(TU,(x^1,\ldots, x^n,y^1, \ldots, y^n))$, and its dual Euclidean metric $\left< \cdot, \cdot \right>_{\ast}$ on $(T^{\ast}U, (x^1, \ldots, x^n, \alpha_1, \ldots, \alpha_n))$. 
  Denote their corresponding $C^0$-Finsler structures by $\Vert \cdot \Vert$ and $\Vert \cdot \Vert_{\ast}$.
  Let $\beta(x) = (\beta_1(x^1, \ldots, x^n), \ldots,$ $\beta_n(x^1, \ldots, x^n))$ be the representation of $\beta$ in coordinate systems.
  We have that $\{\beta (x) = F_{\ast}(x_0,\alpha_0)\}$ is a supporting hyperplane of $B_F[x,0,1]$ at $X_{u(v)}(x)=(1,0,$ $\ldots, 0)$, what implies that $\beta_1(x) = F_{\ast}(x_0,\alpha_0)$ and $F_{\ast}(x,\beta(x)) = F_{\ast}(x_0,\alpha_0)$ for every $x$.
  For a sufficiently small neighborhood $V$ of $x_0$, there exist $c > 0$ such that $B_{F}[x,0,1] \supset B_{\Vert \cdot\Vert}\left[x,0,c \right]$ for every $x \in V$, what implies that $B_{\Vert \cdot \Vert_{\ast}}\left[x, 0,\frac{F_{\ast}(x_0,\alpha_0)}{c}\right] \supset B_{F_\ast}[x,0,F_{\ast}(x_0,\alpha_0)]$ for every $x \in V$.
  This means that $\Vert (\beta_1(x), \ldots, \beta_n(x)) \Vert_{\ast} \leq \frac{F_{\ast}(x_0,\alpha_0)}{c}$, that is, this family of functionals indexed by $x \in V$ is bounded.
  The fact that $\beta$ isn't continuous at $x_0$ implies that there exists $\varepsilon > 0$ and a sequence $x_n \rightarrow x_0$ in $V$ such that  $\Vert \beta(x_n) - \beta(x_0)\Vert_{\ast} > \varepsilon$ for every $n\in \mathbb{N}$ (Here we make the canonical identification among cotangent spaces in Euclidean spaces). 
  In addition, the fact that $\beta$ is bounded implies that, eventually taking a subsequence, we have $\beta(x_n) \rightarrow \xi \in T^{\ast}_{x_0}M \backslash 0$, $F_{\ast}(x_0,\xi) = F_{\ast}(x_0,\alpha_0)$ and $\Vert \xi - \alpha_{0}\Vert_{\ast} \geq \varepsilon$.
  Therefore, $\alpha_{0}(X_{u(v)}(x_0)) = \xi (X_{u(v)}(x_0)) = F_{\ast}(x_0,\alpha_0)$ and $F_{\ast}(x_0, \xi) = F_{\ast}(x_0,\alpha_0)$, what implies that $\{\alpha_0 = F_{\ast}(x_0,\alpha_0)\}$ and $\{\xi = F_{\ast}(x_0,\alpha_0)\}$ are distinct supporting hyperplanes of $B_F[x_0,0,1]$ at $v$, what gives a contradiction.
  Therefore, $\beta$ is continuous at $x_0$, what settles the proof.
\end{proof}

\begin{remark}
  \label{remark_C1_boundary}
It is a classical result from convex analysis that a nonempty, compact and convex subset of $\mathbb{R}^n$ with nonempty interior has a $C^1$-hypersurface as its boundary if and only if every boundary point admits a unique supporting hyperplane (see \cite{Schneider}, Theorem 2.2.4). 
\end{remark}

\begin{corollary}
  \label{corollary_unique_hyperplane}
  Let $(M,F)$ be a Pontryagin-type $C^0$-Finsler manifold such that $S_F[x,0,1]$ is a $C^1$-hypersurface of $T_xM$ for every $x \in M$, that is, $F_x: T_{x}M\backslash\{0\} \rightarrow \mathbb{R}$ is a $C^1$ map. 
  Then $\mathcal{E}$ depends only on $F$, and every minimizing path of $M$ parametrized by arc length is an extremal.
\end{corollary}

\begin{proof}
  $\mathcal{E}$ depends only on $F$ due to Proposition \ref{proposition_unique_hyperplane} and Remark \ref{remark_C1_boundary}.
  
  Let us prove the second statement.
  We know that $F_{\ast}(x(t),\alpha(t))$ is constant for a Pontryagin extremal $(x(t),\alpha(t))$ due to PMP. 
  Moreover, for $c>0$, $(x(t),\alpha(t))$ is a Pontryagin extremal if and only if $(x(t),c \, \alpha(t))$ is also a Pontryagin extremal due to \eqref{align_extended_gedesic_field}, \eqref{align_differential_inclusion} and \eqref{align_Hamiltonian_vector_field_coordinates}.  Therefore, there isn't any loss of generality if we suppose that $F_{\ast}(x(t),\alpha(t))=1$, because the existence of such an extremal is equivalent to the existence of a Pontryagin extremal $(x(t),c.\alpha(t))$, for any $c>0$.  With this assumption, there exist a bijection between $\alpha \in S_{F_{\ast}}[x,0,1]$ and supporting hyperplanes $\{ \alpha=1 \}$ of $B_F[x,0,1]$.

  Let $x(t)$ be a minimizing path in $M$. 
  It is not difficult to prove that when $x(t)$ is restricted to an open subset $U$ where an $F$-family is defined, there exist a unique Pontryagin extremal $(\tilde{x}(t),\tilde{\alpha}(t))$ such that $\tilde{x}(t) = x(t)\vert_{U}$. 
  In fact, if $(\tilde{x}(t),\bar{\alpha}(t))$ is another Pontryagin extremal such that $\tilde{\alpha}(t_0) \neq \bar{\alpha}(t_0)$, then the continuity of $\tilde{\alpha}$ and $\bar{\alpha}$ implies $\tilde{\alpha}(t)\neq \bar{\alpha}(t)$ on a neighborhood of $t_0$, what contradicts the uniqueness of the supporting hyperplane $\{\alpha(t)=1\}$ of $\tilde{x}^{\prime}(t)$.
  Therefore, every minimizing path $x(t)$ can be uniquely lifted locally to a Pontryagin extremal $(\tilde{x}(t),\tilde{\alpha}(t))$ such that $F_{\ast}(\tilde{x}(t),\tilde{\alpha}(t)) = 1$.
  The global lift $(x(t),\alpha(t))$ is obtained glueing the local lifts, what is possible due to the uniqueness of the local lifts, what gives a Pontryagin extremal corresponding to $x(t)$, and settles this corollary.
\end{proof}

\begin{remark}
  In modern language, the geometric structure that Riemann defined in his seminal work \cite{Riemann} is a function $F:TM \rightarrow \mathbb{R}$ such that $F\vert_{\slittangentbundle}$ is smooth and $F\vert_{T_xM}$ is a norm for every $x\in M$ (see also \cite{ChernNotices} and \cite{Spivak2}).
  These structures are of Pontryagin type due to their smoothness. 
  Later, Finsler geometry was developed under a different hypothesis, where $F\vert_{T_{x}M}$ is an asymmetric norm such that the Hessian of $F^{2}\vert_{T_xM\backslash \{0\}}$ at $y$ is positive definite for every $y \in T_xM\backslash\{0\}$ (see \cite{Finsler}). 
  Corollary \ref{corollary_unique_hyperplane} establishes that $\mathcal{E}$ generalizes the geodesic spray to a broader class of structures, including the Riemann's original one.
\end{remark}

\bibliographystyle{amsplain}

\bibliography{\myreferences}{}

\end{document}